\documentclass[12pt]{article}
\usepackage{amssymb}

\newtheorem{thm}     {Theorem}[section]
\newtheorem{prop}    [thm]{Proposition}

\newtheorem{lemma}   [thm]{Lemma}

\newcommand{\proof} {\noindent{\bf Proof. }}

\newcommand{\B}{\mathbb B}
\newcommand{\C}{\mathbb C}
\newcommand{\D}{\mathbb D}

\newcommand{\R}{\mathbb R}

\newcommand{\st}{{\rm st}}

\def\Re{{\rm Re\,}}

\def\bar{\overline}

\begin{document}

\title{Regularized maximum of
strictly plurisubharmonic functions on an almost complex manifold  }
\author{Alexandre Sukhov{*} }
\date{}
\maketitle

{\small

* Universit\'e des Sciences et Technologies de Lille, Laboratoire
Paul Painlev\'e,
U.F.R. de
Math\'e-matique, 59655 Villeneuve d'Ascq, Cedex, France, sukhov@math.univ-lille1.fr
The author is partially suported by Labex CEMPI.

}
\bigskip

{\small Abstract. We prove that the maximum of two smooth strictly plurisubharmonic functions on an almost complex manifold can be uniformly approximated by smooth strictly plurisubharmonic functions.



\bigskip

\section{Introduction}

In this  note the following result is proved:

\begin{prop}
\label{theo1}
 Let $u_1$ and $u_2$ be smooth strictly $J$-plurisubharmonic functions on an  almost complex manifold $(M,J)$. Then for every $\varepsilon > 0$ and every relatively compact domain $\Omega \subset M$ there exists a function $\tilde u$ smooth and strictly $J$-plurisubharmonic in $\Omega$  such that 
\begin{eqnarray}
\label{est}
\max\{ u_1, u_2 \} \leq \tilde u \leq \max\{ u_1, u_2 \} + \varepsilon
\end{eqnarray} 
 on $\Omega$.
\end{prop}

This property of plurisubharmonic functions is standard in the case where an almost complex structure $J$ is integrable.  In their recent lecture notes Cieliebak and Eliashberg \cite{CE} raised this question in the almost complex case. Proposition  \ref{theo1} gives an affirmative answer. It is a corollary of a more precise result  (Theorem \ref{smoothing}). The elementary proof   is proceeded  by the known construction of the regularized maximum function (see the monograph of Demailly \cite{De}) and requires only minor modifications with respect to the integrable case.

\section{Plurisubharmonic functions on almost complex manifolds: the background}

For convenience of reader's I recall  related  notions.

{\bf 2.1. Almost complex structures.} Consider a smooth (everywhere this means of class $C^\infty$) manifold $M$ of dimension $2n$. An almost complex structure $J$ is a smooth map   assigning to each point $p \in M$ a linear isomorphism $J(p): T_pM \to T_pM$ of the tangent space $T_pM$ such that $J(p)^2 = -I$; here   $I:T_pM \to T_pM$ denotes the identity map. A couple $(M,J)$ is called an almost complex manifold of complex dimension n. The present paper concerns with  almost complex manifolds and structures  of class $C^\infty$ though the results still hold under a lower regularity.

A $C^1$-map $f:M' \to M$ between two almost complex manifolds  $(M',J')$ and $(M,J)$, is called  $(J',J)$-complex or  $(J',J)$-holomorphic  if it satisfies {\it the Cauchy-Riemann equations} $df \circ J' = J \circ df$. By the elliptic regularity, such a map is necessarily of class $C^\infty$. If the source manifold  $M'$ is a Riemann surface, the
holomorphic maps are called $J$-complex (or $J$-holomorphic) curves.
Denote by  $\D$   the
unit disc in $\C$ and  by $J_{\st}$   the standard complex structure
of $\C^n$; the value of $n$ will be clear from the context. In the case where $M' = \D$ and   $J' = J_{st}$,  a $J$-holomorphic map $f$ is called a $J$-{\it complex  disc}.

\bigskip

{\bf 2.2. Local coordinates.} Let $\B_n$ denotes the Euclidean unit ball of $\C^n$.  For every point $p$ in an almost complex manifold  $M$ of complex dimension $n$, there exist a neighborhood $U$ of $p$ and a
coordinate diffeomorphism $z: U \to \B_n$ with  $z(p) = 0$, such that  the
direct image $z_*(J):= dz \circ J \circ dz^{-1}$ satisfies  $z_*(J)(0) = J_{st}$. Choose an integer  $k \geq 1$  and $\lambda_0 > 0$. Composing $z$ with isotropic dilations in $\C^n$ one can   additionly achieve the condition
$\vert\vert z_*(J) - J_{st}\vert\vert_{C^k(\overline {\B}_n)} \leq \lambda_0$.

In these local coordinates  
$J$ is represented by a $\R$-linear operator
$J(z):\C^n\to\C^n$, $z\in \C^n$ such that $J(z)^2=-I$. We  use the notation $\zeta=\xi+i\eta\in\D$ for the standard complex coordinate in $\C$. Then the Cauchy-Riemann equations for a $J$-complex  disc $f:\D\to\B_n$ have the form $\partial_\eta f=J(f)\partial_\xi f$. Similarly to  \cite{Aud},
the structure $J$ defines a  unique smooth  complex $n\times n$
matrix function $A_J=A_J(z)$ allowing to write  the Cauchy-Riemann equations as an elliptic quasilinear deformation of the usual $\overline\partial$-equation:
\begin{eqnarray}
\label{CR}
\partial_{\bar\zeta} f + A_J(f)\partial_{\bar\zeta} \bar f = 0
\end{eqnarray}
 $A_J$ is called the {\it complex matrix} of $J$ in the coordinates $z$. If $z':U \to \C^n$ is another coordinate chart and $A_J'$ is the complex matrix of $J$ in the coordinates $z'$, then $A_J' = ((\partial_{z} z')A_J - \partial_{\overline z} z')(\partial_{\overline z}\overline{z}' - (\partial_z\overline{z}') A_J)^{-1}$
(see, for example, \cite{DiSu}).  The condition $J(0) = J_{st}$ means that
\begin{eqnarray}
\label{CR1}
A_J(0) = 0
\end{eqnarray}
As it was mentioned, given $k \geq 0$ the norm $\parallel A_J \parallel_{C^k(\B_n)}$ can be made arbitrarily small by an isotropic dilation of coordinates.

\bigskip

{\bf 2.3. Plurisubharmonic functions.} Let  $u$ be a real $C^2$ function on an open subset $\Omega$ of an almost complex manifold  $(M,J)$. Denote by $J^*du$ the
 differential form acting on a vector field $X$ by $J^*du(X):= du(JX)$. Given point $p \in M$ and a tangent vector $V \in T_p(M)$ consider  a smooth vector field $X$ in a
neighborhood of $p$ satisfying $X(p) = V$. 
The value of the {\it complex Hessian} ( or  the  Levi form )   of $u$ with respect to $J$ at $p$ and $V$ is defined by $H(u)(p,V):= -(dJ^* du)_p(X,JX)$.  This definition is independent of
the choice of a vector field $X$. For instance, if $J = J_{st}$ in $\C$, then
$-dJ^*du = \Delta u d\xi \wedge d\eta$; here $\Delta$ denotes the Laplacian. In
particular, $H_{J_{st}}(u)(0,\frac{\partial}{\partial \xi}) = \Delta u(0)$.

Recall some  basic properties of the complex Hessian (see for instance, \cite{DiSu}):

\begin{lemma}
\label{pro1}
Consider  a real function $u$  of class $C^2$ in a neighborhood of a point $p \in M$.
\begin{itemize}
\item[(i)] Let $F: (M',J') \longrightarrow (M, J)$ be a $(J',
  J)$-holomorphic map, $F(p') = p$. For each vector $V' \in T_{p'}(M')$ we have
$H_{J'}(u \circ F)(p',V') = H_{J}(u)(p,dF(p)(V'))$.
\item[(ii)] If $f:\D \longrightarrow M$ is a $J$-complex disc satisfying
  $f(0) = p$, and $df(0)(\frac{\partial}{\partial \xi}) = V \in T_p(M)$ , then $H_J(u)(p,V) = \Delta (u \circ f) (0)$.
\end{itemize}
\end{lemma}
 Property (i) expresses the holomorphic invariance of the complex Hessian. Property (ii) is often useful in order to compute the complex Hessian  on a given 
tangent vector $V$.

 Let $\Omega$ be a domain $M$. An upper semicontinuous function $u: \Omega \to [-\infty,+\infty[$ on $(M,J)$ is
{\it $J$-plurisubharmonic} if for every $J$-complex disc $f:\D \to \Omega$ the composition $u \circ f$ is a subharmonic function on $\D$. By Proposition
\ref{pro1}, a $C^2$ function $u$ is plurisubharmonic on $\Omega$ if and only if it has    a 
positive semi-definite complex Hessian on $\Omega$ i.e. $H_J(u)(p,V) \geq 0$  for any $ p \in \Omega $ and $V \in T_p(M)$. The equivalence of these two definitions still holds in the general case if the complex Hessian is understood in the sense of currents. This was established by Pali \cite{Pa} for continuous functions and by Harvey-Lawson \cite{HL} for  upper semicontinuous functions.

A real $C^2$ function $u:\Omega \to \R$ is called {\it strictly $J$-plurisubharmonic} on $\Omega$, if $H_J(u)(p,V) > 0$ for each $p \in M$ and $V \in T_p(M) \backslash \{ 0\}$. Obviously, these notions   are local: an upper semicontinuous (resp. of class $C^2$) function on $\Omega$ is  $J$-plurisubharmonic (resp. strictly) on $\Omega$ if and only if it is $J$-plurisubharmonic  (resp. strictly) in some open neighborhood of each point of $\Omega$.

As above, choosing local coordinates near $p$ we may identify a neighborhood
of $p$ with a neighborhood of the origin, assuming that $J$-complex discs
are solutions of the equations (\ref{CR}) and the condition (\ref{CR1})  holds. The following technical result shows that after an additional change of local coordinates one can achieve a further normalization of the  complex matrix $A_J$ of $J$.

\begin{lemma}
\label{normalization}
There exists   local coordinate  diffeomorphism fixing the origin whose components are polynomials of degree at most 2 and   with the linear part equal to the identity map, such that in the new coordinates $z$ the following holds:
\begin{itemize}
\item[(i)]  the matrix function $A_J$  from the equation (\ref{CR}) still satisfies (\ref{CR1}) and additionly satisfies the condition 
\begin{eqnarray}
\label{norm}
\partial_z A_J(0)  = 0
\end{eqnarray}
\item[(ii)] Every real function $u = u(z)$ of class $C^2$  satisfies
\begin{eqnarray}
\label{Hessian}
H_J(u)(0;V) = H_{J_{st}}(u)(0;V) = \sum_{k,j=1}^n \frac{\partial^2 u}{\partial z_k \partial \overline{z}_j}(0)V_k\overline{V}_j
\end{eqnarray}
 for each vector $V = (V_1,...,V_n) \in \C^n$.
\end{itemize}
\end{lemma}

I learned this result from unplublished lecture notes of Chirka although it is possible that it was known before. An elementary proof can be found, for instance, in \cite{DiSu}. 

Local coordinates given by Lemma \ref{normalization} near a point $p \in M$ are called {\it adapted coordinates} at $p$. The existence of adapted local  coordinates provides a convenient way to evaluate the complex Hessian of a $C^2$ function at a given point of an almost complex manifold. 

Conclude this section by two  remarks.

1. The condition (\ref{CR1}) alone (without (\ref{norm}) in general does not guarantee the property (\ref{Hessian}). Consider, for example, the harmonic function $u(\zeta) = \Re \zeta$ in a neighborhood of the origin in $\C$. After the  local change of coordinates $\zeta = F(\tilde\zeta) = \tilde\zeta + \vert \tilde \zeta \vert^2$ we obtain the function $\tilde u(\tilde\zeta) = u(F(\tilde\zeta)) =
\Re\tilde\zeta + \vert \tilde \zeta \vert^2$. Since $dF(0) = I$, the complex structure $J = (F)^*(J_{st}):=(F^{-1})_*(J_{st})$ satisfies the condition (\ref{CR1}). However $\tilde u$ is not strictly $J$-plurisubharmonic since $H_J(\tilde u)$ vanishes identically in a neighborhood of the origin.

2. Denote by $\parallel \bullet \parallel$ the Euclidean norm on $\C^n$. In the adapted coordinates the function $u(z) = \parallel z \parallel^2$ is strictly $J$-plurisubharmonic and satisfies (\ref{Hessian}).  If local coordinates are not adapted and only the condition (\ref{CR1}) holds, then  in general the equality  (\ref{Hessian}) can fail.  However, even in this case $u$ is still strictly $J$-plurisubharmonic near the origin if  the norm $\parallel A \parallel_{C^1}$ is small enough.

\section{Regularized max-function} In this section Proposition \ref{theo1} is proved. I literally follow the presentation given by Demailly \cite{De}.

Denote by  $d\mu^k$  the standard Lebesgue measure on $\R^k$. Fix a $C^\infty$ function $\omega:\R \longrightarrow \R^+$  such that the support of $\omega$ is contained in $[0,1]$ and $\int_\R \omega d\mu^1 = 1$.

Given $\theta = (\theta_1,...,\theta_k) \in (\R^+)^k$ consider  {\it the regularized max-function}: 

\begin{eqnarray}
\label{Max}
M_\theta(t_1,...,t_k) = \int_{\R^k} \max\{t_1 + s_1,...,t_k+s_k\} \prod_{j=1}^k \theta_j^{-1}\omega(s_j/\theta_j)d\mu^k(s)
\end{eqnarray}
defined for $t = (t_1,...,t_k) \in \R^k$. Its utility arises from the following

\begin{lemma}
\label{MaxLemma}
For every $\theta \in (\R^+)^k$ we have:
\begin{itemize}
\item[(i)] The function $t=(t_1,...,t_k) \mapsto M_\theta(t_1,...,t_k)$  is smooth on $\R^k$.
\item[(ii)] $\max\{t_1,...,t_k\} \leq M_\theta(t_1,...,t_k) \leq \max \{t_1 + \theta_1,...,t_k + \theta_k\}$ for every $t \in \R^k$.
\item[(iii)] $M_\theta(t_1+a,...,t_k +a) = M_\theta(t_1,...,t_k) + a$ for every $a \in \R$ and 
every  $t \in \R^k$.
\item[(iv)] Let  $u_1,...,u_k$ be smooth  $J$-plurisubharmonic functions on $M$.  Then $M_\theta(u_1,...,u_k)$ is a smooth $J$-plurisubharmonic function on $M$. 
\end{itemize}
\end{lemma}
\proof (i) Perfoming the change of variables $t_j + s_j = s_j'$ we obtain
$$M_\theta(t_1,...,t_k) = \int_{\R^k} \max\{s_1',...,s_k'\} \prod_{j=1}^k \theta_j^{-1}\omega((s_j'-t_j)/\theta_j)d\mu^k(s')$$
Differentiation of this integral with respect to the parameters $t_j$ shows that $M_\theta$ is smooth. (ii) and (iii) follow from the assumptions on the function $\omega$. At (iv), consider a $J$-complex disc $f:\D \to M$. For all $s_j$ the composition 
$$\D \ni \zeta \mapsto \max\{u_1(f(\zeta)) + s_1,...,u_k(f(\zeta))+s_k\}$$
is a subharmonic function on $\D$. Therefore the function 
$$\D \ni \zeta \mapsto M_\theta(u_1(\zeta),...,u_k(\zeta))$$
is subharmonic in $\D$  by Theorem 2.4.8 in \cite{R}. Hence the function 
$M_\theta(u_1,...,u_k)$ is a smooth $J$-plurisubharmonic function on $\Omega$. $\blacksquare$

Choose a continuous hermitian metric $h$ on $(M,J)$ and denote by $h_p(X)$ the value of  $h$ at point $p \in M$ on a vector $V \in T_pM$. Proposition \ref{theo1} is a consequence of the following more precise result.

\begin{thm}
\label{smoothing}
Let $\alpha$ be a smooth real function on $M$.  Suppose that $u_1,...,u_k$ are smooth functions on $\Omega$ satisfying $H_J(u_j)(p,V) \geq \alpha (p)h_p(V)$ for every point $p \in \Omega$ and each vector $V \in T_pM$. Then for every $\theta = (\theta_1,...,\theta_k) \in (\R^+)^k$ the function $\tilde u= M_\theta(u_1,...,u_k)$ is smooth and satisfies  $H_J(\tilde u)(p,V) \geq \alpha(p) h_p(V)$ on $\Omega$.
\end{thm}
\proof  Fix a point $p$ and $\delta > 0$. Choose adapted coordinates at the point $p$. In particular,  in these coordinates $p$ corresponds to the origin. Let us simply write $J$ instead of $(z)_*(J)$. Since the coordinates are adapted, the function $z \mapsto h_0(z)$ is a positive definite $J_{st}$-hermitian form and $H_J(h_0)(0,V) = H_{J_{st}}(h_0)(0,V) = h_0(V)$.
Consider the smooth functions  $v_j(z) = u_j(z) - \alpha(0)h_0(z) + \delta \parallel z \parallel^2$. They  satisfy $H_{J}(v_j)(0,V) = H_{J_{st}}(v_j)(0,V) \geq \delta \parallel V \parallel^2$. Hence they are $J$-plurisubharmonic in a neighborhood of the origin ( in general,  depending on $\delta$). By (iv)  of  Lemma \ref{MaxLemma} the smooth function
$\tilde v:= M_\theta(v_1,...,v_k)$  is $J$-plurisubharmonic near the origin as well and by (iii) Lemma \ref{MaxLemma}  one has $\tilde v = \tilde u - \alpha(0)h_0(z) + \delta \parallel z \parallel^2$.
Since $\delta > 0$ is arbitrary, we obtain that  $H_J(\tilde u)(0,V) \geq \alpha(0)h_0(V)$.  $\blacksquare$

In order to prove Propposition \ref{theo1} it suffices now for a given $\varepsilon > 0$ to choose $\theta_j = \varepsilon$, $ j=1,2$ and to apply Theorem \ref{smoothing} with $\alpha$ equal to the constant function $1$. The property  (ii) of Lemma \ref{MaxLemma} implies  the estimate (\ref{est}). The proof is completed.


{\footnotesize

\begin{thebibliography}{CIT}




\bibitem{Aud} M.Audin, J.Lafontaine (Eds.) "Holomorphic curves in symplectic geometry",
Birkhauser, Progress in Mathematics, V.117 (1994).

\bibitem{CE} K.Cieliebak, Y.Eliashberg, {\it Stein structures: existence and flexibility}. 
Lecture notes, available at {\tt http://www.renyi.hu/$\sim$cast2012/stein.pdf}


\bibitem{De} J.-P. Demailly,  "Complex analytic and differential geometry", Electronic publication available at
{\tt http://www-fourier.ujf-grenoble.fr/$\sim$demailly/books.html}. 


\bibitem{DiSu} K.Diederich, A.Sukhov, {\it Plurisubharmonic exhaustion functions and almost complex Stein structures}, Mich. Math. J. {\bf 56}(2008), 331-355.

\bibitem{HL} F.R. Harvey, H.B.Lawson, Jr., {\it Potential theory on almost complex manifolds},
arXiv: 1107.2584v2 [math. CV].


\bibitem{Pa} N. Pali, {\it Fonctions plurisousharmoniques et courants positifs de type (1,1) sur une vari\'et\'e presque complexe}, Manuscripta math. {\bf 118} (2005), 311-337.

\bibitem{R} Th. Ransford, {\it Potential theory on the complex plane}, Cambridge Univ. Press, 1995.






\end{thebibliography}
}
\end{document}